\documentclass[a4paper, 
11pt
]{article}

\usepackage{subfig}
\usepackage{pgfplots}
\usepackage{pgfplotstable}
\usepackage{mathtools}
\usepackage{multicol}
\usepackage{comment}
\usepackage{booktabs}
\pgfplotsset{compat=1.5}
\usepackage{amssymb}
\usepackage{url}
\usepackage{bm}

\usepackage{pdfsync}
\usepackage{float}
\usepackage{tabularx}
\usepackage{enumerate}
\usepackage{array}
\usepackage{xspace}

\usepackage{authblk}

\begin{document}

\title{Shape Optimization by means of Proper Orthogonal Decomposition and Dynamic Mode Decomposition}

\author[1,2]{Nicola~Demo\footnote{nicola.demo@sissa.it}}
\author[1]{Marco~Tezzele\footnote{marco.tezzele@sissa.it}}
\author[2]{Gianluca~Gustin\footnote{gianluca.gustin@fincantieri.it}}
\author[2]{Gianpiero~Lavini\footnote{gianpiero.lavini@fincantieri.it}}
\author[1]{Gianluigi~Rozza\footnote{gianluigi.rozza@sissa.it}}

\affil[1]{Mathematics Area, mathLab, SISSA, International School of Advanced Studies, via Bonomea 265, I-34136 Trieste, Italy}
\affil[2]{Fincantieri - Divisione Navi Mercantili e Passeggeri, Cantieri Navali Italiani SpA, Trieste, Italy}

\maketitle

\begin{abstract}
Shape optimization is a challenging task in many engineering fields,
since the numerical solutions of parametric system may be
computationally expensive. This work presents a novel optimization
procedure based on reduced order modeling, applied to a naval hull
design problem. The advantage introduced by this method is that the
solution for a specific parameter can be expressed as the combination
of few numerical solutions computed at properly chosen parametric
points. The reduced model is built using the proper orthogonal
decomposition with interpolation (PODI) method.
We use the free form deformation (FFD) for an automated perturbation
of the shape, and the finite volume method to simulate the multiphase
incompressible flow around the deformed hulls. Further computational
reduction is done by the dynamic mode decomposition (DMD) technique:
from few high dimensional snapshots, the system evolution is
reconstructed and the final state of the simulation is faithfully
approximated. Finally the global optimization algorithm iterates over
the reduced space: the approximated drag and lift coefficients are
projected to the hull surface, hence the resistance is evaluated for
the new hulls until the convergence to the optimal shape is achieved.
We will present the results obtained applying the described procedure
to a typical Fincantieri cruise ship.
\end{abstract}

\section{Introduction}
\label{sec:intro}
Thanks to the improved capabilities in terms of computational
infrastructures in the last decade, the classical simulation-based
design (SBD) has evolved into an automatic optimization procedure,
that is simulation-based design optimization (SBDO). For what concerns
shape optimization, the SBDO procedure consists of three main
features: a shape parametrization and deformation tool, an high
fidelity solver, and an optimization algorithm.
In this work we present a complete SBDO pipeline of the bulbous bow of
a hull advancing in calm
water. It consists in an efficient shape parametrization technique, an
high fidelity solver based on the finite volume method, two different
model reduction techniques in order to speed up both the single high
fidelity simulation and the optimization procedure, and
an optimization algorithm. For what concerns shape
optimization of hulls we cite among
others~\cite{d2017nonlinear,diez2015design}.

In the framework of reduced order modelling (ROM) and efficient
shape parametrization techniques, a first possible choice is related to the shape morphing
method itself. Since there are no particular constraints, we focus on
the so-called \emph{general purpose} methods, and among them we mention \emph{Free
  Form Deformation} (FFD)~\cite{sederbergparry1986,LassilaRozza2010},
\emph{Radial Basis Functions} (RBF) interpolation~\cite{buhmann2003radial,morris2008cfd,manzoni2012model,tezzele2017combined} or
\emph{Inverse Distance Weighting} (IDW) interpolation~\cite{forti2014efficient, BallarinDAmarioPerottoRozza2017}. Generally
speaking, these methods involve the displacement of some
control points in order to induce a deformation on the domain, and we identify
the parameters as the displacements of the control points. 

As high fidelity solver we use a finite volume method based one. We
remark that this particular choice does not represent a constraint,
since the pipeline proposed is independent and the solver can be
considered as a black box. We simply need the flow fields of each
simulation at particular time intervals. Each simulation is
accelerated by the dynamic mode decomposition (DMD)
technique~\cite{schmid2010dynamic,schmid2011applications,kutz2016dynamic,le2017higher}. 
Originally introduced in the fluid mechanics community, the DMD has
emerged as a powerful tool for analyzing the dynamics of nonlinear
systems. We use it to get an estimate of the total resistance at
regime simulating actually only few seconds.
Finally a reduced space constructed with the POD modes is employed by
the optimization algorithm for a fast evaluation of the total
resistance for new parameters. In particular we use an interpolation
based approach (POD with interpolation) where the new solution is
obtained by interpolating the low rank solutions into the parametric
space. This non intrusive choice allows us to apply the pipeline to
different shape optimization problems, changing only the high fidelity
solver or the parametrization technique. 
In this work we present the results of the optimization of the bulbous bow of
a cruise ship designed and built by Fincantieri.

\section{The Wave Resistance Approximation of an Hull Advancing in
  Calm Water}

Let $\mathbb{D} \subset \mathbb{R}^m$ with $m \in \mathbb{N}$, be the
set of parameters, and assume that it is a box in $\mathbb{R}^m$. 
Let $\mathcal{M}(\boldsymbol{x}; \boldsymbol{\mu}): \mathbb{R}^n \to
\mathbb{R}^n$ be a shape morphing function, mapping the reference
domain $\Omega$ to the deformed domain $\Omega(\boldsymbol{\mu})$,
that is $\Omega(\boldsymbol{\mu}) = \mathcal{M}(\Omega; \boldsymbol{\mu})$.
In the following of this work (see Section~\ref{sec:ffd}) we present
the actual definition of $m$ and $\mathcal{M}$ for this problem. 

The result of the fluid dynamic simulations will heavily depend on the specific
shape of the hull. In this work we will assess the effect of the shape
deformation on the total resistance, that is the main fluid dynamic performance
parameter. In particular for each point in the parameter space $\mathbb{D}$, a
morphing of the original hull geometry is created thanks to the function
$\mathcal{M}$. The high fidelity solver based on finite volumes (see
        Section~\ref{sec:solver}) takes these geometries as input and performs
the simulation of the flow past the hull advancing in calm water for each
specific hull. To come up with a resistance estimate at regime we adopt a
reduction strategy that is based on the dynamic mode decomposition presented in
the Section~\ref{sec:dmd}. With all the snapshots collected in this part of the
pipeline, we construct a reduced space through a proper orthogonal
decomposition with interpolation method. This space is used to perform fast
evaluations of the total resistance for new parameters by the optimization
algorithm.

\section{Shape Morphing Based on Free Form Deformation}
\label{sec:ffd}
Free form deformation (FFD) is a widely used deformation
technique both in academia and in industry. In this section we
summarize the FFD-based morphing technique. For a further insight on
the original formulation see~\cite{sederbergparry1986}, and for more
recent works the reader can refer
to~\cite{lombardi2012numerical,rozza2013free,sieger2015shape,salmoiraghi2016isogeometric,forti2014efficient,salmoiraghi2017,tezzele2017dimension}.
All the algorithms have been implemented in the open source Python package
PyGeM~\cite{pygem}, which is used to perform the shape morphing in the
numerical results showed in Section~\ref{sec:results}. 

Basically the idea of FFD is to embed the part of the geometry we want to morph
in a lattice and to deform it using a trivariate
tensor-product of B\'ezier or B-spline functions. Thus, by moving only
the control points of such lattice, we can produce a continuous and
smooth deformation. The FFD procedure
can be subdivided into three steps. First, we need to map the physical
domain $\Omega$ to the reference domain $\widehat{\Omega}$ through the
map~$\psi$. Then, we move some control points
$\boldsymbol{P}$ of the lattice to achieve the desired deformation,
using the map $\widehat{T}$. The displacement of such points are the
weights of the FFD and they represent the parameters
$\boldsymbol{\mu}$. Finally we apply the back mapping from the deformed
reference domain to the deformed physical domain $\Omega(\boldsymbol{\mu})$ by the map $\psi^{-1}$.
So we can express the FFD map $T$ by the composition of these three
maps, that is
$T(\cdot, \boldsymbol{\mu}) = (\psi^{-1} \circ \widehat{T} \circ \psi)
(\cdot, \boldsymbol{\mu})$.

\begin{figure}[thb]
\centering
\includegraphics[width=.48\textwidth, trim=0 155 0 0]{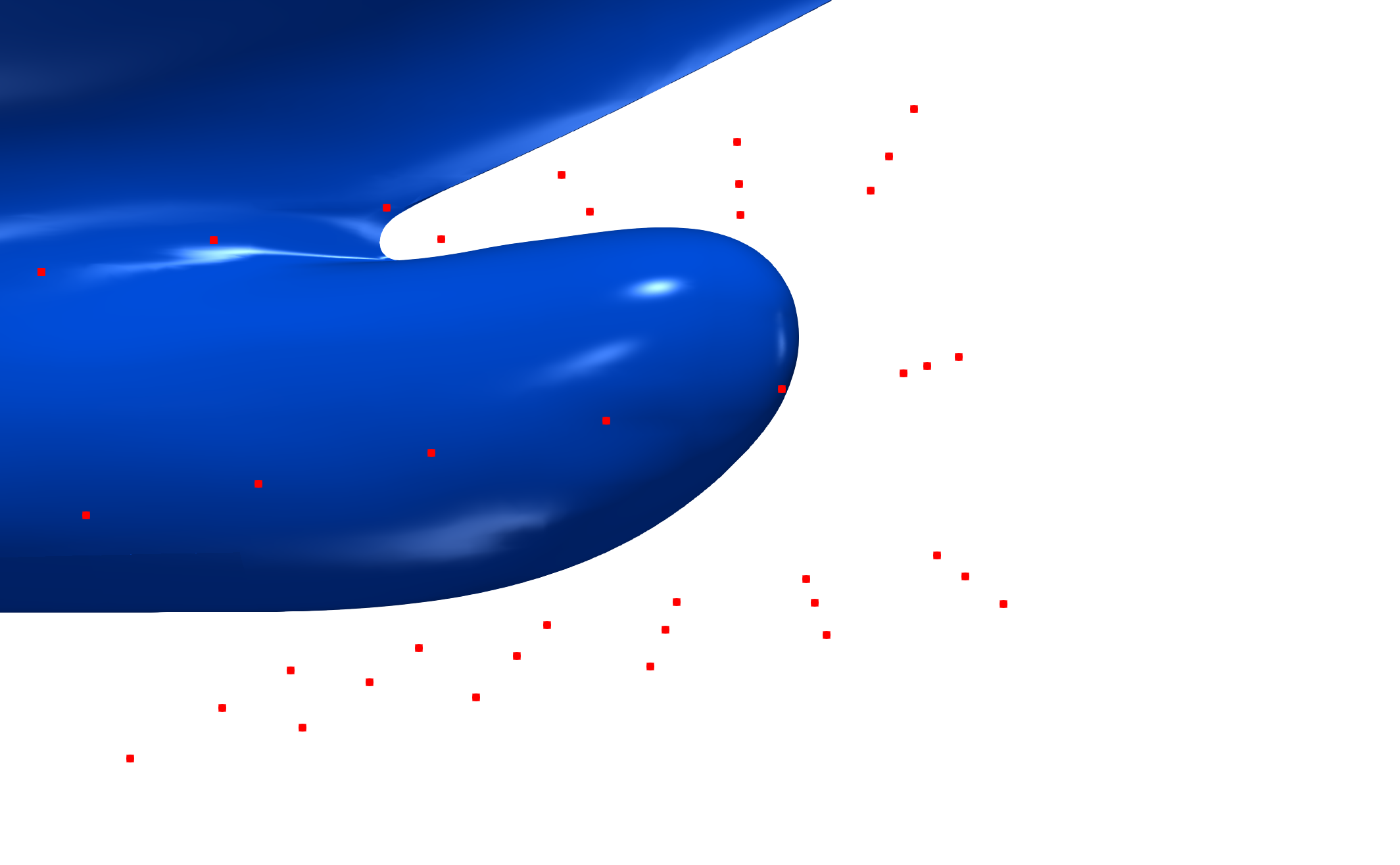}
\includegraphics[width=.48\textwidth, trim=0 155 0 0]{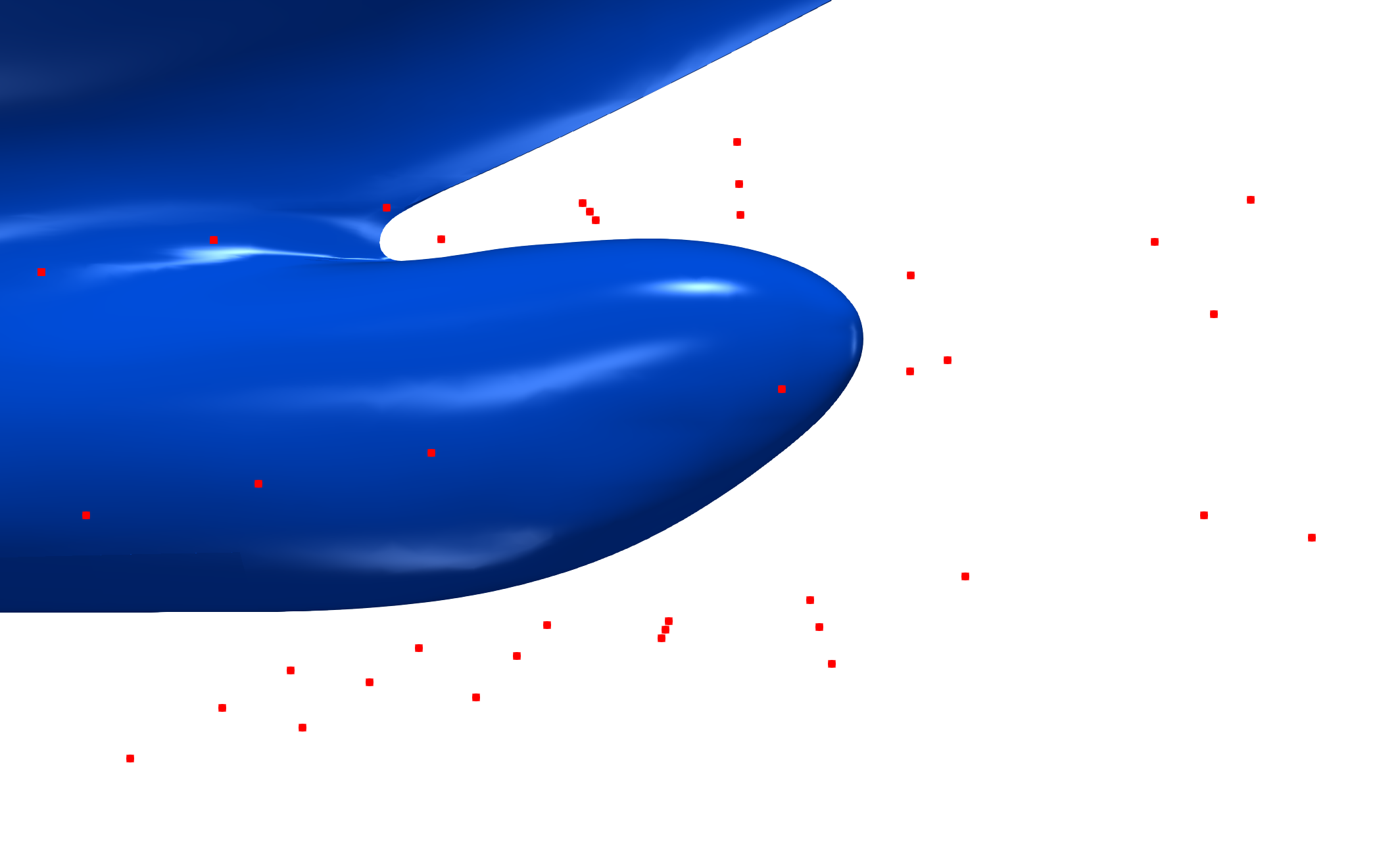}
    \caption{Bulbous bow and FFD points, before (left) and after (right) a deformation.}
\label{fig:deform}
\end{figure}
In Figure~\ref{fig:deform} it is possible to see the original bulbous
bow on the left and an example of a deformed one on the right,
respectively with the original FFD control points and with the
stretched ones.

\section{High Fidelity Solver Based on Finite Volume Method}
\label{sec:solver}
Dealing with turbulent flows, the physical system is described through
the Reynolds-averaged Navier Stokes (RANS) equations. In order to approximate
the turbulent fluid, these equations decompose the instantaneous velocity into
fluctuating and time-averaged parts \cite{reynolds}. We denote respectively
these quantities as $u'(x, t)$ and $\bar{u}(x, t)$, such that $u(x, t) = u'(x,
        t) + \bar{u}(x, t)$. Let us put this decomposition to the incompressible
continuity equation and to the momentum equation to obtain the RANS equations:
\begin{equation}
\begin{dcases}
    \frac{\partial{\bar{u}_i}}{\partial{x_i}} = 0\\[-8pt]
\frac{\partial{\bar{u}_i}}{\partial{t}} + 
    \bar{u}_j\frac{\partial{\bar{u}_i}}{\partial{x_j}} =
    -\frac{1}{\rho}\frac{\partial{\bar{p}}}{\partial{\bar{x}_i}} +
    \nu \frac{\partial^2{\bar{u}_i}}{\partial{x_j}\partial{x_j}} -
    \frac{\partial{\overline{u'_iu'_j}}}{\partial{x_j}}
\end{dcases}
\end{equation}
The introduction of Reynold stress term $\overline{u'_iu'_j}$ requires an
additional modelling in order to close the system solve it: in our work
we use the popular $SST~k-\omega$ model \cite{menter1994two}.

In the present work, the governing equations are solved using the finite volume
method (FVM). This technique is widely spread in the computational fluid
dynamics community since it guarantees the conservation of all the quantities
being based on the conservative form of the equations. Moreover, it is easily
applicable on complex unstructured meshes. Basically, to provide the numerical
solution of the partial differential equations, the space domain is subdivided
into a finite number of non-overlapping polyhedra called finite volumes. The
governing equations are so integrated on each finite volume and the integral
values are approximated on the reference cells. For the numerical approximations
of the turbulent flow over the hull surface, we use the C++ FVM-based library
OpenFOAM~\cite{of}.

\section{Spatial and Temporal Reduction Using Dynamic Mode Decomposition}
\label{sec:dmd}

Dynamic mode decomposition (DMD) is inspired by and closely related to
Koopman-operator analysis~\cite{rowley2009spectral}. DMD was developed by Schmid
in~\cite{schmid2010dynamic}, and since then has emerged as a powerful tool for
analyzing the dynamics of nonlinear systems, and for postprocessing
spatio-temporal data in fluid
mechanics~\cite{schmid2011application,schmid2011applications,stegeman2015proper}.
DMD has gained popularity in the fluids community, primarily because it
provides information about the dynamics of a flow, and it is applicable
even when those dynamics are nonlinear~\cite{rowley2009spectral}. Many variants
of the DMD arose in the last years like multiresolution DMD, forward backward
DMD, compressed DMD, and higher order DMD. For a complete review refer
to~\cite{kutz2016dynamic,le2017higher}.  Since DMD relies only on the
high-fidelity measurements, like experimental data and numerical simulations,
    it is an equation-free algorithm, and it does not make any assumptions
    about the underlying system.

Let us consider a sequential set of data vectors $\{\boldsymbol{x}_1, \dots, \boldsymbol{x}_m\}$
where $\boldsymbol{x}_k \in \mathbb{R}^n$ for all $k \in [1, m]$. We assume that the vectors are
sampled from a continuous evolution $\boldsymbol{x} (t)$, and
equispaced in time. Hence, the $\boldsymbol{x}_k$
represents the state of the system. We also suppose that
the dimension $n$ of a snapshot is larger than
the number of snapshots $m$, that is $n>m$. The basic idea is that
there exists a linear operator $\mathbf{A}$ (also called Koopman
operator) that approximates the
nonlinear dynamics of $\boldsymbol{x} (t)$, that is $\boldsymbol{x}_{k+1} = \mathbf{A} \boldsymbol{x}_k$.
The DMD modes and eigenvalues are intended to approximate the eigenvectors and
eigenvalues of $\mathbf{A}$.  To obtain the minimum approximation error across
all these \textit{snapshots}, it is possible to arrange them in two matrices
such that 
$\mathbf{X} = \begin{bmatrix} \mathbf{x}_1 & \dotsc & \mathbf{x}_{m-1}\end{bmatrix}$
and
$\mathbf{Y} = \begin{bmatrix} \mathbf{x}_2 & \dotsc & \mathbf{x}_{m} \end{bmatrix}$,
with
$\mathbf{x}_i =\begin{bmatrix} x_i^1 & \dotsc & x_i^n \end{bmatrix}^\intercal$.
We underline that each column of $\mathbf{Y}$ contains the state
vector at the next timestep of the one in the corresponding
$\mathbf{X}$ column. We want to find $\mathbf{A}$ such that the relation between the
matrices $\mathbf{X}$ and $\mathbf{Y}$ is $\mathbf{Y} \approx \mathbf{A} \mathbf{X}$.
The best-fit $\mathbf{A}$ matrix is given by $\mathbf{A} = \mathbf{Y} \mathbf{X}^\dagger$,
where the symbol $^\dagger$ denotes the Moore-Penrose pseudo-inverse. 
Since $n>m$, the matrix $\mathbf{A}$ has
$n^2$ elements and it is difficult to decompose it or to
handle it. The DMD algorithm projects the data onto a low-rank subspace
defined by the POD modes --- the first left-singular vectors of the matrix
$\mathbf{X}$ computed by the truncated SVD, that is $\mathbf{X} \approx
\mathbf{U}_r \bm{\Sigma}_r \mathbf{V}^*_r$. We call $\mathbf{U}_r$
the unitary matrix whose columns contain the first $r$ modes. The low-dimensional operator is built as:
$
\mathbf{\tilde{A}} = \mathbf{U}_r^* \mathbf{A} \mathbf{U}_r =
\mathbf{U}_r^* \mathbf{Y} \mathbf{V}_r \bm{\Sigma}_r^{-1}
$.
So we can avoid the explicit calculation of the high-dimensional
operator~$\mathbf{A}$, obtaining the matrix $\mathbf{\tilde{A}} \in
\mathbb{C}^{r\times r}$.
We can now reconstruct the eigenvectors and eigenvalues of the matrix
$\mathbf{A}$ thanks to the eigendecomposition of $\mathbf{\tilde{A}}$
as $\mathbf{\tilde{A}} \mathbf{W} = \mathbf{W} \bm{\Lambda}$.
In particular the elements in $\bm{\Lambda}$ correspond to the nonzero
eigenvalues of $\mathbf{A}$, while the real eigenvectors, the so called
\textit{exact} modes~\cite{tu2014dynamic}, can be computed as $\bm{\Phi} =
\mathbf{Y}\mathbf{V}_r \bm{\Sigma}_r^{-1} \mathbf{W}$. 
This algorithm has been implemented `in house', as well as many of its variants,
in an open source Python package called PyDMD~\cite{demo18pydmd}.

\section{Proper Orthogonal Decomposition with Interpolation}
Reduced order methods (ROMs) have become a fundamental tool for the complex
systems analysis. They make possible a remarkable reduction of the
computational cost in the calculation of a solution of a parametric partial
differential equation (PDE). In this section, we discuss the model
reduction based on the proper orthogonal decomposition (POD), focusing on the
non intrusive approach.

The basic idea is to represent the system as a linear combination of few main
structures, the so called \textit{modes}. Using the POD method, these structures
are the orthogonal basis functions individuated through the modal analysis on
the solutions vectors~\cite{hesthaven2016certified}. Hence, the high-fidelity
solutions have to be computed, for some proper parameters, by the full-order
solver, using the desired accuracy and the desired number of degrees of freedom.
This first phase is the most computationally expensive, due to the calculation
and the storage of several high-order solutions (\textit{offline phase}). Then,
in the \textit{online phase}, these solutions are combined for an efficient and
reliable approximation of the solution for a new generic parameter. In the POD
Galerkin strategy --- intrusive method --- the discretized PDEs are projected
on the space spanned by the POD modes and the significantly smaller system
obtained is solved. For further details about this methodology and examples of
its application, we cite~\cite{BallarinRozza2016, Stabile2017, Stabile2018}.
Conversely, in a POD interpolation (PODI) approach --- non intrusive method ---
the solutions are projected on the low dimensional space spanned by the POD
modes. In this way, the solutions are described as linear combination of the
modes: the coefficients of this combination, the so called \textit{modal
coefficients}, are interpolated in order to provide the coefficients for any
point belonging to the parameter space. With the interpolated coefficients and
the modes we can finally compute the approximated solution. PODI method has the
advantage to require only the
solutions vectors, with no assumption on the underlying system, though the
accuracy of the approximation depends on the interpolation method chosen.
\begin{figure}[t]
\centering
\includegraphics[width=.77\textwidth, trim=0 40 0 0]{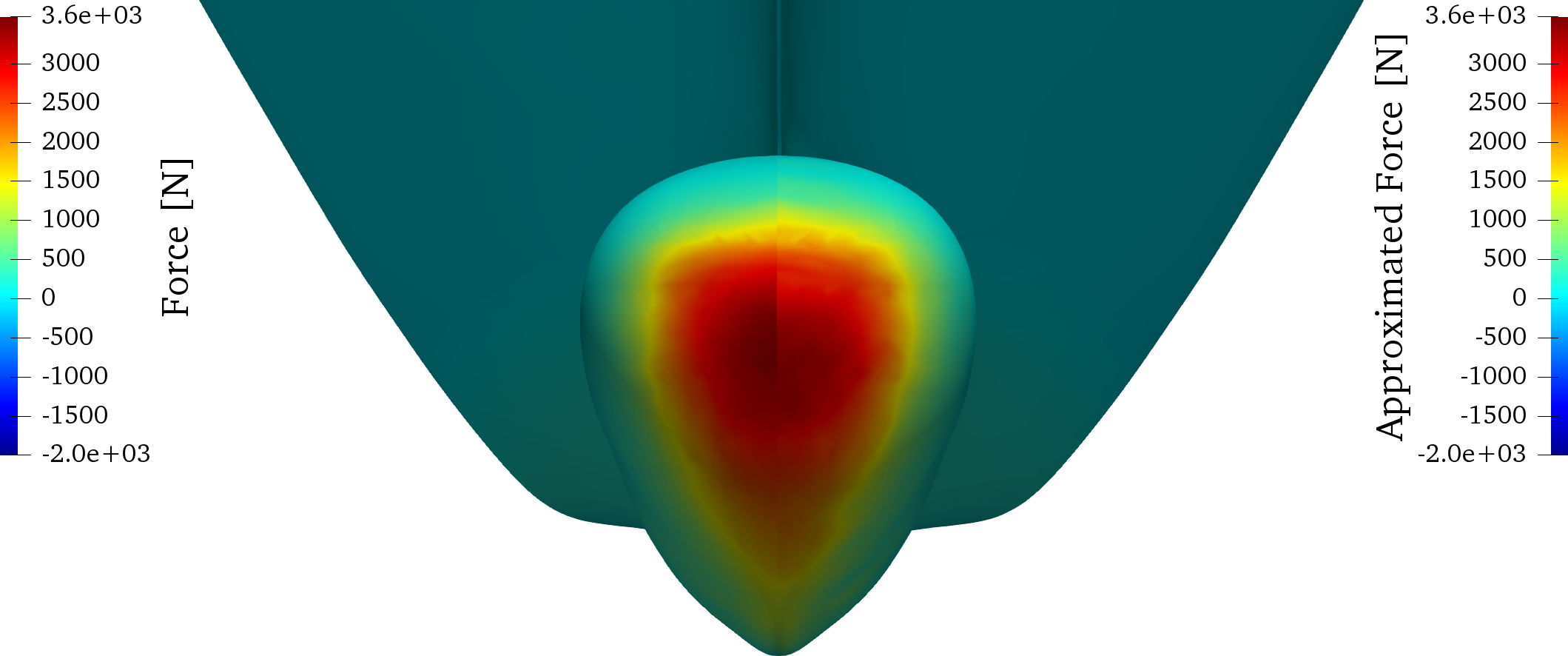}
\caption{A comparison between the numerical solution from the full-order model (left) and the numerical solution from the reduced model using a PODI strategy. The picture represents the force on the hull domain.}
\label{fig:bulbo-rom}
\end{figure}

Within this work, we adopted the PODI method in order to maintain the
optimization pipeline completely independent from the full-order model. The
method has been implemented in the open source Python package
EZyRB~\cite{demo18ezyrb}. In Figure~\ref{fig:bulbo-rom} we show an example of the
numerical solution computed with this software, compared with the high-fidelity
validation from the full-order model.

\section{Numerical Results}
\label{sec:results}
\begin{figure}[t]
\frame{\includegraphics[width=.30\textwidth, trim=0 20 0 0]{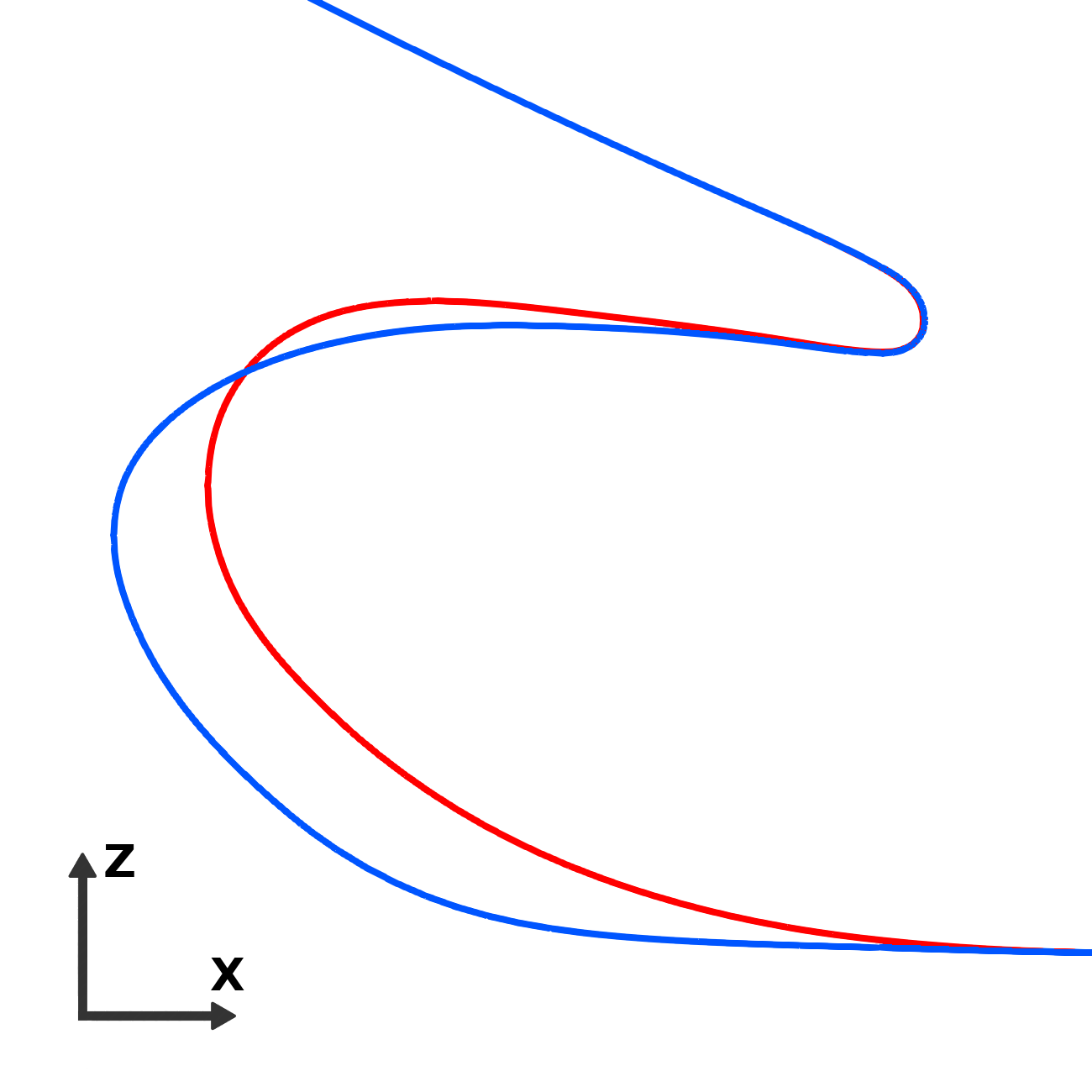}}
\hfill
\frame{\includegraphics[width=.30\textwidth, trim=0 20 0 0]{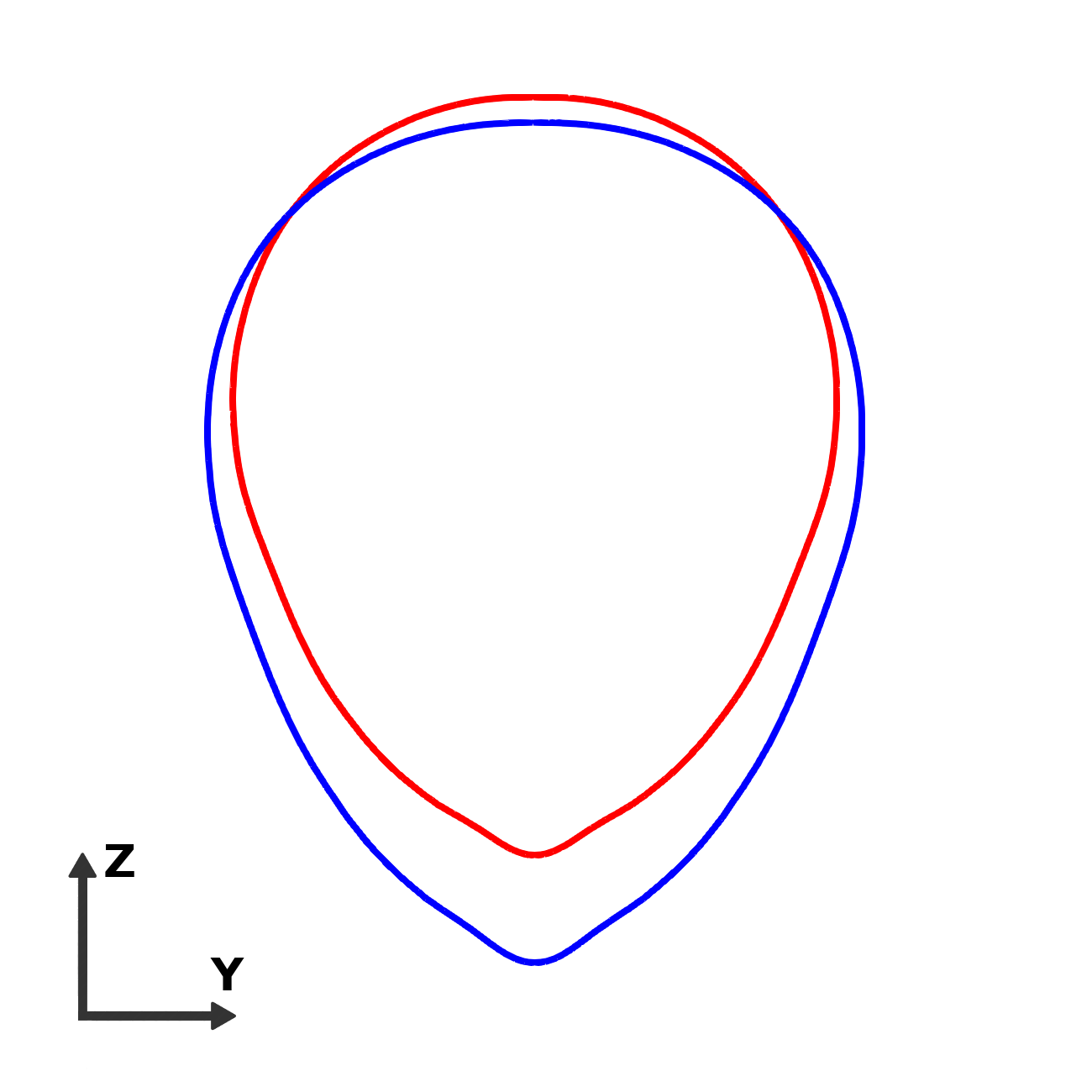}}
\hfill
\frame{\includegraphics[width=.30\textwidth, trim=0 20 0 0]{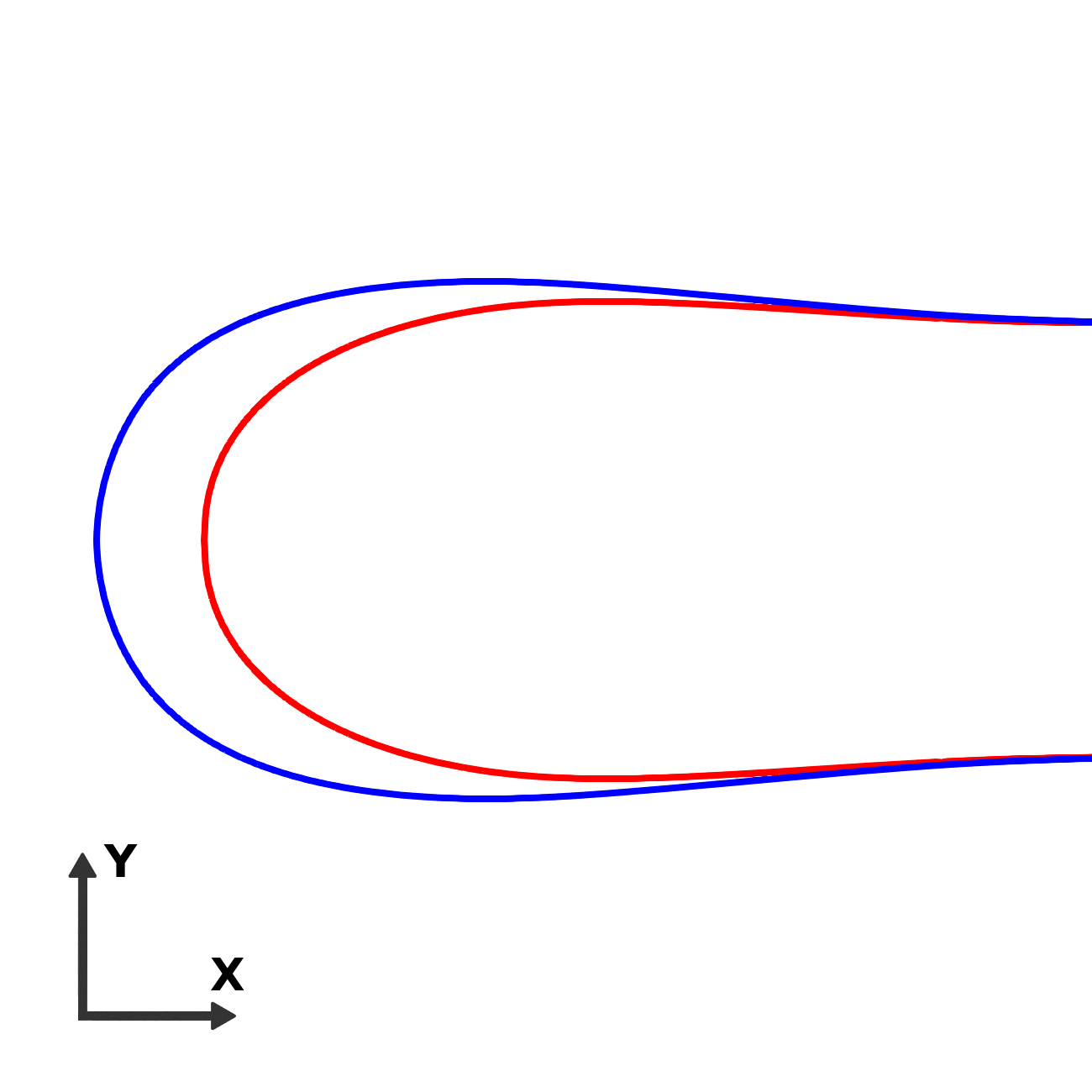}}
\caption{Three sections of the bulbous bow with the original shape (red) against the optimized one (blue).}
\label{fig:bulbopt}
\end{figure}

In this section we present the numerical results obtained by the application of
the above described pipeline on a Fincantieri cruise ship. We adopted the FFD
technique for the surface deformation: we used 5 different parameters and we
set the parameter space dimensions to obtain physically meaningful
hulls. We generated 62 different deformed shapes, 32 corresponding to the
vertices of the parameter space and 30 corresponding to uniform sampling
points in the parameter space. We simulated the advance in calm water of all
these deformed hulls at a constant speed corresponding to Fr = 0.2 by using the FV method, collecting 20 snapshots of the
full-order simulation between the 50th and 60th second. These snapshots,
containing the pressure and shear stress fields on the hull surface, have
been used to estimate the hull resistance at regime by the DMD
algorithm. Exploiting the correlation between the parametric points and
the reconstructed solutions, we created the reduced space with the PODI strategy, interpolating the modal coefficients with \textit{radial basis functions}
with \textit{multiquadric} kernel. The low
computational cost of a single parametric solution allows us to apply
expensive global optimization algorithms, such as \textit{surrogate-based
optimization}. This algorithm generates several designs of experiment
and interpolates the objective function evaluated in these points to create
a surrogate model, then this model is evaluated until convergence to the optimal
point. The procedure iterates until the accuracy of the surrogate
model is reached. This method allows a quasi real-time optimization.
Figure~\ref{fig:bulbopt} the optimized bulb compared to the original one. After the
optimization loop, the best hull has been used as input for a further
high-fidelity simulation, to validate the numerical solution of the
reduced model. The automatic optimization procedure reached a
remarkable reduction in the hull resistance: comparing the full-order
solutions, the optimized ship resistance is 2\% lower compared to the original
ship, while the error between the high-fidelity and the PODI solution,
 evaluated in the optimal parametric point, is around the 8\%. Despite
the error is bigger than the resistance reduction, we underline that we mainly focused on the minimization of the total
resistance, hence the error was not involved in the stop criteria. We recall that
the accuracy of the reduced model depends on the richness of the database of
the high-fidelity solutions. Thus the PODI method can reach the wanted
precision enriching the database, for example by adding iteratively the
high-fidelity validation computed in the optimal point. 
Regarding this work, the reached reduction was enough to demonstrate a concrete
application of the optimization pipeline.

\section{Conclusions and Perspectives}
\label{sec:the_end}

In this work we presented a shape optimization pipeline combining DMD
and POD with interpolation. We applied it to the problem of minimizing
the total resistance of a hull advancing in calm water varying the
shape of the bulbous bow of a Fincantieri cruise ship. The numerical
results show that with only 5 parameters we achieved a reduction of the resistance of 2\%, remaining independent from the solver. Several enhancements
on the pipeline could be foreseen. Among the possible ones, we mention
a deeper investigation on the interpolation method and on the
optimization algorithm itself.

\section*{Acknowledgements}
This work was partially funded by the project HEaD, ``Higher Education
and Development'', supported by Regione FVG --- European Social Fund FSE
2014-2020, and by European Union
Funding for Research and Innovation --- Horizon 2020 Program --- in
the framework of European Research Council Executive Agency: H2020 ERC
CoG 2015 AROMA-CFD project 681447 ``Advanced Reduced Order Methods
with Applications in Computational Fluid Dynamics'' P.I. Gianluigi
Rozza.

\end{document}